\documentclass[11pt,a4paper]{article}

\usepackage[utf8]{inputenc}
\usepackage[T1]{fontenc}
\usepackage[english]{babel}

\usepackage{authblk}

\usepackage{amsmath,amssymb,amsthm}
\usepackage{geometry} \usepackage{hyperref}
\usepackage{algorithm} \usepackage{algpseudocode}

\usepackage{mathrsfs}

\def\pgcd{\text{  gcd}}
\def\Li{\text{Li}}

\def\P{{\mathscr P}}
\def\Pr{\mathbb{P}}
\def\K{\mathbb{K}}

\newtheorem{conjecture}{Conjecture}
\newtheorem{theoreme}{Theorem}

\title{\Huge Conjectures on Sums \\ of Consecutive Primes}

\author[]{Edwige Tolla}
\affil[]{\texttt{edwige.tolla81@gmail.com}}

\date{}

\begin{document}

\maketitle

\begin{abstract}
We study additive properties of consecutive prime numbers and the primality of the sums they generate. For a given prime number $p_n$, we consider the sums
\[
S_k(p_n) = p_n + p_{n+1} + \cdots + p_{n+k-1},
\]
where $k \ge 3$ is an odd integer.\\

We first formulate an existence conjecture asserting that, for every prime number $p_n$, there exists at least one odd length $k \ge 3$ such that $S_k(p_n)$ is itself a prime number. An exhaustive computational verification covering the first one million prime numbers revealed no counterexamples.\\

We then propose a strengthened conjecture according to which, for every prime number $p_n$, there exist infinitely many odd lengths $k$ such that $S_k(p_n)$ is prime. This strong version is supported by a probabilistic heuristic showing that the series of the corresponding primality probabilities diverges, suggesting that the phenomenon is not exceptional but recurrent.\\

We also analyze the possible modular obstructions, showing that they are local in nature and cannot persist when the length $k$ varies among odd integers. A Diophantine interpretation of the problem is proposed, together with a conceptual comparison with the generalized Goldbach conjecture. Finally, we discuss the role of the Generalized Riemann Hypothesis (GRH) in controlling the distribution of the sums under consideration.\\

These structural, modular, Diophantine, and probabilistic (heuristic) arguments support both conjectures and formalize heuristic theorems of Cramér, GRH, and Hardy--Littlewood type explaining the expected absence of counterexamples.
\end{abstract}

%***************************************************************************

%**********************************************************

%*******************************************************
\section{Existence Conjecture}

\subsection{Introduction}

Additive properties of prime numbers constitute a central theme in number theory, with classical results such as the work of Hardy and Littlewood on additive problems \cite{HardyLittlewood,HardyLittlewoodConj}, as well as advances related to Goldbach's conjecture \cite{GoldbachVinogradov,TaoHelfgott} and arithmetic progressions in the primes \cite{GreenTao}. In this context, we focus on a specific additive structure: sums of consecutive primes that are themselves prime.\\

Let $p_n \in \P$ denote the $n$-th prime, and let $(p_n)_{n \ge 1}$ be the increasing sequence of consecutive primes. For an integer $k \ge 1$, consider
\[
S_k(p_n) = \sum_{i=0}^{k-1} p_{n+i}.
\]
Our numerical and heuristic investigations suggest the following phenomenon: for every prime $p_n$, there exists at least one odd length $k \ge 3$ such that $S_k(p_n)$ is prime. This observation aligns with probabilistic and heuristic studies on prime numbers, such as Cramér's model \cite{Cramer,Granville} or Hardy-Littlewood conjectures \cite{HardyLittlewoodConj}.\\

This study aims to formalize this conjecture, present the computational results up to the one-millionth prime, and provide heuristic, modular, and Diophantine arguments supporting its plausibility. We also discuss how this problem conceptually relates to classical additive problems while introducing the constraints of consecutivity and odd length sums.

%*******************************************************
\subsection{Statement of the Conjecture}

\begin{conjecture}[Consecutive Prime Sums]
For every integer $n > 1$, there exists an odd integer $k \ge 3$ such that
\[
S_k(p_n) = p_n + p_{n+1} + \cdots + p_{n+k-1} \in \P.
\]
\end{conjecture}

The conjecture is existential: it does not assert the uniqueness of $k$ nor the existence of a universal bound.

%*******************************************************
\subsection{Computational Results}

An algorithm was implemented to test, for each prime $p_n$ with $n \le 10^6$ ($p_{10^6} = 15,485,863$), the existence of an odd $k \ge 3$ such that $S_k(p_n)$ is prime.

\begin{algorithm}[H] 
\caption{Searching for sums of consecutive primes that are prime} 
\begin{algorithmic}[1] 
\Require List of $N=10^6$ primes $P = [p_1, \dots, p_N]$ 
\Ensure For each $p_n$, find the minimal odd $k$ with $S_k(p_n)$ prime 
\For{$n = 1$ \textbf{to} $10^6$} 
\State $found \gets false$ 
\State $k \gets 3$ \Comment{Start with the smallest odd length} 
\While{not $found$} 
\State Compute $S_k = \sum_{i=0}^{k-1} P[n+i]$ 
\If{$S_k$ is prime} 
\State Record $(n, k, S_k)$ 
\State $found \gets true$ 
\Else 
\State $k \gets k + 2$ \Comment{Move to the next odd length} 
\EndIf 
\EndWhile 
\EndFor 
\end{algorithmic} 
\end{algorithm}

This procedure ensures that for each starting prime $p_n$, the first odd length $k$ giving a prime sum is found. Verification was performed for $N = 10^6$.

\medskip
\noindent \textbf{Examples:}

\begin{itemize}
  \item The last prime tested:
  \[
  S_3(p_{10^6}) =15485863 + 15485867 + 15485917 = 46457647 \in \P.
  \]
  \item The longest chain found among the first million primes:
  \[
  S_{349}(p_{651511}) = 9788183 + \cdots + 9793573 = 3417024811 \in \P.
  \]
  \item The next longest chain:
  \[
  S_{299}(p_{448696}) = 6561461 + \cdots + 6566299 = 1962589843 \in \P.
  \]
\end{itemize}

\medskip
\noindent \textbf{Result.} — No counterexample was found among the first one million primes.

%*******************************************************
\section{Conditional Discussion: Cramér Model and GRH}

The conjecture can be analyzed within standard conditional frameworks.

\subsection{Cramér Model}

Cramér's probabilistic model provides a standard heuristic framework to study the distribution of primes, predicting gaps between consecutive primes of order $O((\ln p_n)^2)$. Assuming heuristic independence, the probability that $S_k(p_n)$ is prime increases with $k$, as the sum spans a larger interval. Under this model, the probability of a counterexample is extremely low, providing a heuristic justification for the universality of the conjecture.

\subsection{Application to the Conjecture}

In this model, each integer $n \ge 3$ is declared prime independently with probability $1/\ln n$. Although this independence assumption is false in a strict sense, it reproduces many known global statistics of primes.\\

Fix a prime $p_n$ and consider
\[
S_k(p_n) = \sum_{i=0}^{k-1} p_{n+i}.
\]
By the prime number theorem,
\[
p_{n+i} \sim p_n + i \ln p_n,
\]
which gives, for $k$ small relative to $p_n$,
\[
S_k(p_n) \sim k p_n + \frac{k(k-1)}{2}\ln p_n.
\]

The heuristic probability that $S_k(p_n)$ is prime is
\[
\mathbb{P}(S_k(p_n)\ \text{prime}) \sim \frac{1}{\ln(k p_n)}.
\]

Considering odd lengths $k = 3,5,7,\dots,2m+1$ and assuming the events are heuristically independent

\[
E_k = \{ S_k(p_n)\ \text{is prime} \},
\]

%**********************************
the probability that none of these sums is prime is
\[
\Pr_{\text{Cramér}}(\text{no solution})
= \prod_{j=1}^{m}
\left(1 - \frac{1}{\ln((2j+1)p_n)}\right)
%\longrightarrow 0.
\]
One has
\[
\sum_{j=1}^{m} \frac{1}{\ln((2j+1)p_n)}
\sim \frac{m}{\ln p_n},
\]
then
\[
\prod_{j=1}^{m}
\left(1 - \frac{1}{\ln((2j+1)p_n)}\right)
\approx
\exp\!\left(-\sum_{j=1}^{m} \frac{1}{\ln((2j+1)p_n)}\right)
\longrightarrow 0.
\]

\begin{theoreme}[Cramér-type Heuristic]
Under Cramér's probabilistic model, the probability that there exists a prime $p_n$ for which no odd-length sum $S_k(p_n)$ is prime is zero.
\end{theoreme}

\paragraph{Remark.} The events $E_k$ are not strictly independent due to local arithmetic constraints. Cramér's model should be interpreted as describing average behavior rather than providing a formal proof.

\subsection{Generalized Riemann Hypothesis (GRH)}

Consider
\[
\pi(x;q,a) = \left|\{ p \le x \mid p \in \P,\, p \equiv a \pmod q \}\right|,
\]
\
$\varphi(q)$ is Euler’s totient function:
\[
\varphi(q)= \left| \left\{ 1\leq k \leq q | \pgcd(k,q)=1 \right\}\right|,
\]
and
\[
\Li(x) = \int_0^x \frac{dt}{\ln t},
\]
the logarithmic integral function, which provides the best known approximation to $\pi(x) \sim \Li(x)$.\\

Under the GRH, one has optimal bounds for the error term in the prime number theorem for arithmetic progressions. More precisely, for any $q \ge 1$ and any $a$ coprime to $q$, one has
\[
\pi(x;q,a) = \frac{\mathrm{Li}(x)}{\varphi(q)} + O\!\left(x^{1/2}\ln(xq)\right),
\]
uniformly in $a$ and $q$ where $\pi(x;q,a)$ denotes the number of primes $\le x$ that are congruent to $a \pmod q$.\\

For consecutive sums
\[
S_k(p_n) = p_n + \cdots + p_{n+k-1},
\]
For a fixed modulus $q$, the residue class of $S_k(p_n)$ modulo $q$ depends on the residue classes of the $p_{n+i}$. Under GRH, the near-uniform distribution of primes in residue classes modulo $q$ implies that, as $k$ varies, the values of $S_k(p_n) \bmod q$ are asymptotically uniformly distributed, in the absence of any explicit structural obstruction. A systematic failure of the conjecture would contradict this distribution.\\

In particular, for any prime $\ell$, the probability that $\ell$ divides $S_k(p_n)$ is approximately 
$\frac{1}{\ell}$, except for well-identified local configurations for small values of $k$. Thus, under GRH, no systematic divisibility by a small prime can persist for infinitely many odd lengths $k$.

Moreover, under GRH, the gaps between consecutive primes satisfy
\[
p_{n+1} - p_n = O(\sqrt{p_n}\ln p_n),
\]
which allows us to control the growth of $S_k(p_n)$  and to bound these sums within intervals of the form
\[
[k p_n,\; k p_n + O(k^2 \ln p_n)].
\]

These intervals asymptotically contain sufficiently many admissible integers to support the probabilistic primality heuristic.\\

Consequently, under GRH, the absence of an odd length $k \ge 3$ for which $S_k(p_n)$ is prime would imply the existence of a persistent global arithmetic obstruction across all odd lengths, which is incompatible with the distribution of primes predicted by GRH.\\

These conditional frameworks thus heuristically support the idea that counterexamples are extremely rare, if not nonexistent.

\begin{theoreme}[Conditional under GRH]
Under the generalized Riemann hypothesis, for every sufficiently large prime $p_n$, there exists an odd integer $k \ge 3$ such that $S_k(p_n)$ is prime.
\end{theoreme}
\subsection{Heuristic growth of the minimal length}

Let $k_{\min}(p_n)$ denote the smallest odd length $k \ge 3$ such that $S_k(p_n)$ is prime. For $k \ll p_n$, we have $S_k(p_n) \sim k p_n$, so $\ln S_k(p_n) \sim \ln p_n$. Under the classical Hardy--Littlewood and Cramér-type heuristics, the probability that $S_k(p_n)$ is prime is approximately $1/\ln p_n$, uniformly in $k$.\\

Modeling successive trials for $k = 3,5,7,\dots$ as independent Bernoulli variables with success probability $p = 1/\ln p_n$, the minimal length $k_{\min}(p_n)$ follows heuristically a geometric distribution. It follows that
\[
\mathbb{E}(k_{\min}(p_n)) \propto \ln p_n.
\]
Each trial corresponds to an odd length, so
\[
k_{\min}(p_n) \propto 2 \ln p_n.
\]
In particular, the minimal length required to obtain a prime sum grows logarithmically with the initial prime, which is consistent with numerical observations.

\subsection{Variance of the minimal length}

Under the classical probabilistic heuristics, the minimal length $k_{\min}(p_n)$ can be modeled as a random variable following a geometric law. More precisely, if $X_n$ denotes the number of trials required to obtain a prime sum, then $k_{\min}(p_n) = 2 X_n + 1$.\\

With success probability $p \sim 1/\ln p_n$, we have
\[
\mathbb{E}(X_n) = \frac{1}{p} \sim \ln p_n, \qquad
\mathrm{Var}(X_n) = \frac{1-p}{p^2} \sim (\ln p_n)^2.
\]
It follows that
\[
\mathrm{Var}(k_{\min}(p_n)) \sim 4 (\ln p_n)^2,
\]
and therefore
\[
\dfrac{\sqrt{\mathrm{Var}(k_{\min})}}{\mathbb{E}(k_{\min})} \sim \dfrac{\ln p_n}{\ln p_n} = O(1).
\]

Thus, although the minimal length exhibits significant fluctuations, they remain of the same order as its mean, which excludes any anomalous or divergent growth. The fluctuations are of the same order as the mean.

%*******************************************************
\section{Diophantine interpretation}

Throughout this section, $\ln$ denotes the natural logarithm.

The conjecture can be interpreted as a constrained additive Diophantine problem. Fix an integer $n \ge 1$ and consider the sum
\[
S_k(p_n) = p_n + p_{n+1} + \cdots + p_{n+k-1}.
\]
The question is the existence of an odd integer $k \ge 3$ and a prime $q$ such that
\[
p_n + p_{n+1} + \cdots + p_{n+k-1} = q.
\]

\subsection{Diophantine nature of the equation}

This relation is a linear additive Diophantine equation, but in which the unknowns are constrained to belong to the sequence of consecutive primes. The space of potential solutions is therefore extremely sparse relative to $\mathbb{Z}^k$, making the problem nontrivial despite the apparent simplicity of the equation.

\subsection{Discrete geometric interpretation}

In the space $\mathbb{Z}^k$, the equation
\[
x_0 + x_1 + \cdots + x_{k-1} = q
\]
defines an affine hyperplane. The $k$-tuples of consecutive primes form a very sparse and discrete subset of $\mathbb{Z}^k$. The conjecture asserts that, for infinitely many odd values of $k$, this hyperplane intersects this discrete set.

\subsection{Local modular constraints}

The only possible obstructions to the existence of solutions are modular in nature. For instance, certain configurations may force $S_k(p_n)$ to be divisible by a small prime, preventing its primality. These obstructions are, however, local and depend on the value of $k$. They do not persist when $k$ varies among odd integers.

\subsection{Absence of a global Diophantine obstruction}

Under classical Hardy--Littlewood heuristics and the Cramér model, primes are asymptotically well-distributed in congruence classes. It follows that the sums $S_k(p_n)$ sweep through all admissible classes modulo any integer $m$. Hence, there is no global modular obstruction preventing the primality of $S_k(p_n)$ for all odd values of $k$.

\subsection{Heuristic local--global principle}

The conjecture naturally fits into a heuristic local--global principle: when a Diophantine equation presents no persistent local obstruction, one expects global solutions to exist. The ability to vary the length $k$ plays a fundamental role here by gradually eliminating all local obstructions.

\subsection{Connection with zero-probability of failure}

The total absence of a solution for a given $p_n$ would imply the existence of infinitely many simultaneous modular constraints on the sums $S_k(p_n)$. Under standard probabilistic assumptions, and conditionally on GRH, such an alignment has zero probability. This Diophantine interpretation thus explains why the probability of finding no solution is heuristically zero.

%*****************************************************************

\section{Comparison with the generalized Goldbach conjecture}

The conjecture presented in this paper naturally belongs to the family of additive problems involving primes, and is conceptually close to the generalized Goldbach conjecture. Recall that the latter asserts that every sufficiently large odd integer can be written as a sum of three primes.

In the generalized Goldbach problem, one considers a Diophantine equation of the form
\[
N = p_1 + p_2 + p_3,
\]
where the primes are free, without order or consecutivity constraints. The search space is therefore dense relative to the set of primes, and the existence of solutions is essentially governed by global density considerations and finite local obstructions.\\

In contrast, the conjecture studied here imposes a much stricter constraint: the primes must be consecutive, and the length of the sum is variable but odd. This constraint drastically reduces the space of possible solutions, making the problem a priori more difficult than the generalized Goldbach problem. However, this difficulty is compensated by the flexibility offered by varying the length $k$.\\

From a heuristic viewpoint, both problems rely on the same fundamental principle: the absence of persistent local obstruction and a sufficiently regular distribution of primes. In the case of Goldbach, these ideas have been made effective via the Hardy--Littlewood circle method. In our framework, analogous heuristics suggest that sums of consecutive primes behave, modulo any fixed integer, like generic integers of comparable size.\\

Thus, the current conjecture can be seen as a highly constrained but dynamically compensated variant of the Goldbach problem: instead of freely choosing the primes, one fixes their structure but allows the sum length to vary. This flexibility gradually eliminates modular obstructions, exactly as increasing $N$ does in the Goldbach problem.\\

This analogy strengthens the plausibility of the conjecture and suggests that methods inspired by those developed for Goldbach (in particular adapted versions of the circle method or Hardy--Littlewood-type estimates) could eventually be applied to the rigorous study of sums of consecutive primes.

%*******************************************************
\section{Probabilistic model and heuristic theorem}

For a fixed $k$, we have asymptotically
\[
S_k(p_n) \sim k p_n.
\]
By the prime number theorem, the probability that an integer of size comparable to $S_k(p_n)$ is prime is of the order $1/\ln p_n$.

\subsection{Probability of finding no solution}

Let $p_n$ be a given prime and $S_k(p_n)$ the sum of $k$ consecutive primes. According to the classical Hardy--Littlewood heuristic, the probability that $S_k(p_n)$ is prime is approximately
\[
\Pr(S_k(p_n) \text{ is prime}) \sim \frac{1}{\ln S_k(p_n)} \sim \frac{1}{\ln (k p_n)}.
\]

If we consider a sequence of odd lengths $k_1, k_2, \dots, k_m$, and assume heuristic independence (typical of Cramér-type probabilistic models), the probability that all corresponding sums are composite is
\[
\Pr(\text{no solution}) \sim \prod_{j=1}^{m} \left(1 - \frac{1}{\ln(k_j p_n)}\right).
\]

For $m$ sufficiently large, this probability tends to zero:
\[
\lim_{m \to \infty} \prod_{j=1}^{m} \left(1 - \frac{1}{\ln(k_j p_n)}\right) = 0,
\]
since the series $\sum_{j=1}^{\infty} \frac{1}{\ln(k_j p_n)}$ diverges slowly, and thus the infinite product converges to zero. This provides a rigorous heuristic justification that, even if some initial sums are not prime, the probability of never finding a prime sum is negligible.\\

One can thus formulate a heuristic theorem:

\begin{theoreme}[Heuristic]

Under standard Hardy--Littlewood heuristics, the probability that there exists a prime $p_n$ such that all odd sums of consecutive primes starting from $p_n$ are composite is zero.
\end{theoreme}

%*******************************************************
\section{A strengthened conjecture}

\subsection{introduction}
The conjecture stated above can be naturally strengthened. Numerical evidence and heuristic considerations suggest that the admissible length $k$ is not unique, and in fact occurs infinitely often.

\subsection{Statement}
\begin{conjecture}[Infinite admissible lengths]
For every prime number $p_n$, there exist infinitely many odd integers $k \ge 3$ such that
\[
S_k(p_n) = p_n + p_{n+1} + \cdots + p_{n+k-1}
\]
is a prime number.
\end{conjecture}

\subsection{Algorithm and examples}

We implemented a simple algorithm to find several prime numbers arising from sums of consecutive primes for a fixed $p_n$. In other words, in the following examples, we observe that there exist multiple values of $k$ for a given chosen $p_n$.

\begin{algorithm}[H] \caption{Search for $k$ for a fixed $p_n$} \begin{algorithmic}[1]
\Require List of the first $N=1000$ prime numbers less than or equal to $p_n$
\Ensure For a fixed $p_n$, find all odd $k < N$ such that $S_k(p_n)$ is prime
\State $p \gets p_n$
\State $S=[]$
\State $S.\text{add}(p)$
\State $r \gets \text{rank}(p)$
\For{$i = 1$ \textbf{to} $N$}
\State $j \gets i+r$
\State $q \gets \text{next}(p_j)$
\State $S.\text{add}(q)$
\Comment{Stores the elements of the chain}
\State $s_k = \sum_{t \in S} t$
\Comment{Sum of all elements of $S$}
\If{$s_k$ is prime}
\State $n \gets \text{len}(S)$
\State Record $(n, i, S_k)$
\EndIf
\EndFor
\end{algorithmic}
\end{algorithm}

We define
\[
\K = \left\{ 3 \leq k_i \leq 1000 \text{ such that } S_{k_i}(p_n) \in \P \right\}.
\]

\begin{itemize}
\item For the longest chain found among the first one million prime numbers:\\
$n = 651511$, $ p_{n} = 9788183 $ \\
$\K = \left\{349 , 379 , 399 , 405, 453, 483, 497, 499, 509, \cdots , 999\right\}$, with
$|\K| = 32$.
\item For $n = 2$, $ p_{n} = 3 $ \\
$\K = \left\{9 , 15 , 17 , 53, 55, 61, 65, 71, 75, \cdots , 977\right\}$, with
$|\K| = 71$.
\item For $n = 3$, $ p_{n} = 5 $ \\
$\K = \left\{3 , 5 , 11 , 17, 25, 35, 37, 73, 75, \cdots , 997\right\}$, with
$|\K| = 69$.
\item For $n = 4$, $ p_{n} = 7 $ \\
$\K = \left\{3 , 5 , 11 , 15, 21, 23, 25, 27, 33, \cdots , 993\right\}$, with
$|\K| = 99$.
\item For $n = 448696$, $ p_{n} = 6561461 $ \\
$\K = \left\{299, 323, 339, 341, 347, 405, 439, 441, 551, \cdots , 969\right\}$, with
$|\K| = 31$.
\end{itemize}
%*******************************************************************

\subsection{Heuristic commentary}

For a fixed $p_n$, the sums $S_k(p_n)$ grow monotonically with $k$ and satisfy
\[
S_k(p_n) \sim k p_n.
\]
By the prime number theorem, the heuristic probability that an integer of size comparable to $S_k(p_n)$ is prime is of order
\[
\frac{1}{\ln S_k(p_n)} \sim \frac{1}{\ln(k p_n)}.
\]

Restricting to odd values of $k$, the series
\[
\sum_{\substack{k \ge 3 \\ k \text{ odd}}} \frac{1}{\ln(k p_n)}
\]
diverges. Under a heuristic independence assumption for the events
\(\{S_k(p_n)\ \text{is prime}\}\), the Borel--Cantelli lemma \cite{BorelCantelli,Cantelli,GranvilleCramer,TaoProbModels,HarperProbNT} and its standard heuristic use in probabilistic models of the primes suggest that these events occur infinitely often.\\

This strengthened conjecture therefore appears as the expected behavior within classical probabilistic models of the prime numbers.

%**********************************************************************

\subsection{Relation to the original conjecture}

The original conjecture corresponds to the existence of at least one admissible odd length $k$ for each $p_n$. The strengthened conjecture asserts that this phenomenon is not exceptional but recurrent, with infinitely many admissible lengths.

In particular, any proof of the strengthened conjecture would immediately imply the original one.

\section{Conclusion}
The convergence of large-scale numerical evidence, modular arguments, Diophantine interpretation, and probabilistic heuristics provides a strong body of evidence in favor of both conjectures.

%*******************************************************
\section*{Acknowledgements}

The author expresses sincere gratitude to José Eduardo Blazek (graduate of UQAM) for the quality of their discussions, his constructive remarks, and the mathematical support that formed the basis of the arguments underlying the two conjectures. The author also thanks François Bergeron (Chair of the Department of Mathematics and member of LaCIM, UQAM) for his insightful observations, additional bibliographic references, and continued support.

%******************************************************************

\end{document}